\documentclass[10pt]{article}
\usepackage[dvips]{graphicx}
\usepackage{amsfonts}
\usepackage{amssymb}
\usepackage{amsmath}
\usepackage[latin1]{inputenc}
\usepackage{color}
\usepackage{diagbox}
\usepackage[left=2.5cm,top=2.5cm,right=2.5cm]{geometry}
\usepackage[ruled,linesnumbered]{algorithm2e}
\usepackage{authblk}

\begin{document}

\title{A New Hybrid Classical-Quantum Algorithm for Continuous\\
Global Optimization Problems}

\author[1]{Pedro C. S. Lara}
\author[1]{Renato Portugal}
\author[2]{Carlile Lavor}
\affil[1]{Laborat\'{o}rio Nacional de Computa\c{c}\~{a}o Cient\'{\i}fica, Petr\'{o}polis, Brazil,

\tt{\small \{pcslara,portugal\}@lncc.br}} %
\affil[2]{IMECC/Unicamp, Campinas, Brazil,

\tt{\small clavor@ime.unicamp.br}}
\date{}
\maketitle

\begin{abstract}
Grover's algorithm can be employed in global optimization methods providing,
in some cases, a quadratic speedup over classical algorithms. This paper
describes a new method for continuous global optimization problems that uses
a classical algorithm for finding a local minimum and Grover's algorithm to
escape from this local minimum. Simulations with testbed functions and
comparisons with algorithms from the literature are presented.
\end{abstract}

%
%


\section{Introduction}

In general, we could say that global optimization is about finding the best
solution for a given problem. Global optimization algorithms play an
important role in many practical problems. A recent review of different
methods in global optimization can be found in \cite{FC09}. 

Recently, some papers addressed the problem of finding the global minimum of
discrete~\cite{BBW05,LK10, LK11} and continuous functions~\cite{PB02, PB05}
using quantum algorithms, where the method used in the discrete case is an
extension of Dürr and H{ø}yer's (DH) algorithm~\cite{DH99}, which in turn is
based on the quantum search algorithm developed by Boyer \textit{et.~al}.
(BBHT)~\cite{BBHT98}. This was made possible after Lov Grover has discovered
the seminal algorithm for searching one item in an unsorted database with $N$
distinct elements~\cite{Gro97}. Grover's algorithm finds the position of the
desired element querying the database $O(\sqrt{N})$ times, which is a
quadratic improvement with respect to the number of times classical
algorithms query the database.

The goal of this work is to present a new method to solve continuous global
optimization problems. It is a hybrid method that uses an efficient
classical algorithm for finding a local minimum and a quantum algorithm to
escape from that, in order to reach the global minimum.

\ The article is organized as follows. In Sec.~\ref{sec:Opt}, we describe
the relationship between quantum algorithms and global optimization
problems, and we review the DH and Baritompa \textit{et.~al}~(BBW)
algorithms. In Sec.~\ref{sec:3}, we describe the new method proposed in this
paper. In Sec.~\ref{sec:results}, we present the simulations and discuss the
results of this work. Finally, in Sec.~\ref{sec:conclusion}, we present our
conclusions.

\section[Quantum Search and Global Optimization Problems]{Quantum Search and
Global Optimization Problems}

\label{sec:Opt}

We start by describing the problem that Grover's algorithm addresses.
Consider a boolean function $h:\{0,\cdots ,N-1\}\rightarrow \{0,1\}$, where $%
N=2^{n}$ and $n$ is some positive integer, such that 
\begin{equation}
h(x)=\left\{ 
\begin{array}{ll}
1, & \hbox{if $x\in M$;} \\ 
0, & \hbox{otherwise,}%
\end{array}%
\right.  \label{ag_f_x}
\end{equation}%
where $M\subseteq \{0,\cdots ,N-1\}$. The goal is to find an element $%
x_{0}\in M$ by querying function $h(x)$ the least number of times. Grover~%
\cite{Gro97} described a quantum algorithm that finds $x_{0}$ when $|M|=1$,
that is, the only element in $M$ is $x_{0}$. Grover's algorithm finds $x_{0}$
with probability greater than or equal to $1-1/N$ by querying $h$ around $%
\frac{\pi }{4}\sqrt{N}$ times. Internally, the algorithm uses an initial
vector in a Hilbert space that undergoes $\frac{\pi }{4}\sqrt{N}$ rotations
of small angles $\theta $, such that $\sin ({\theta }/{2})=1/\sqrt{N}$. The
initial vector is the normalized vector obtained by adding all vectors of
the canonical basis (also known by computational basis) of the $N$%
-dimensional Hilbert space, which yields a vector of equal entries, the
value of which are $1/\sqrt{N}$. After $\frac{\pi }{4}\sqrt{N}$ Grover's
rotations, the final vector has a large overlap with the vector that
represents the solution. This is the guarantee one needs to be sure that a
measurement of the quantum computer in this final state will yield $x_{0}$
with high probability. One rotation in the algorithm is also called a Grover
iteration. In each Grover iteration, function $h$ is queried one time. The
number of times $h$ is queried is the main parameter to measure the
efficiency of hybrid classical-quantum algorithms.

Boyer \textit{et.~al}. generalized Grover's algorithm in two directions \cite%
{BBHT98}. Firstly, they considered the case $|M|>1$ and showed that the
number of rotations required to find one element in $M$ with probability
greater than or equal to $1-1/N$ is 
\begin{equation*}
\frac{\pi }{4}\sqrt{\frac{N}{|M|}}.
\end{equation*}%
The number of Grover iterations decreases when $|M|>1$ because the dimension
of the subspace of the Hilbert space spanned by elements in $M$ is $|M|$.
One obtains a large overlap between the final vector state of the quantum
computer with the solution-subspace with less Grover rotations. That is a
straightforward generalization of Grover's algorithm.

Secondly, the authors addressed the problem of finding one element in $M$
without knowing a priori the number of elements in $M$. The main problem in
this case is to know what is the best number of rotations. If the algorithm
performs too few or too many rotations, the probability to find the correct
results becomes small. Their strategy is to start the algorithm by
performing a small number of Grover rotations followed by a measurement,
which yields an element $x\in \{0,\cdots ,N-1\}$. One checks whether $h(x)=1$%
. If that fails, start over again and increase the number of Grover
rotations. The key point is to determine the increment rate. Formally, the
BBHT algorithm~\cite{BBHT98}, which finds a marked element when the number
of solutions is not known in advance (unknown $|M|$), can be written in a
pseudo-code as Algorithm~\ref{alg:BBHT}. Boyer \textit{et.~al}. proved that
the expected running time of the algorithm is $O(\sqrt{N}/|M|)$ if each
query to function $h$ is evaluated in unit time. They observed that any
value of $\lambda $ in the range $1<\lambda <4/3$ is allowed.

\begin{algorithm}[!htb]
    \Begin {
     Initialize $m=1$ and set $\lambda=8/7$\;
     Choose an integer $j$ uniformly at random such that $0\le j<m$\;
     Apply $j$ Grover's iterations starting from the initial vector\;
     Perform a measurement (let $x$ be the outcome)\;
     If $h(x)=1$, return the result $x$\;
     Otherwise, set $m$ to $\min\{\lambda m, \sqrt{N}\}$ and go to line 3\;
     }
\caption{BBHT Algorithm}
\label{alg:BBHT}
\end{algorithm}

Using the BBHT algorithm, Dürr and H{ø}yer~\cite{DH99} proposed an algorithm
to find the minimum element of a finite list $L$, which can be seen either
as finding the index of the smallest element in a database or as a discrete
global optimization problem. At the beginning, the algorithm selects at
random one element $y$ in $L$ and searches for elements in set $%
M=\{y^{\prime }\in L|y^{\prime }<y\}$ using the BBHT algorithm. If it
succeeds, $y^{\prime }$ is the new candidate for minimum value and the BBHT
algorithm will be used again and again until the minimum is found with
probability greater than 1/2. In order to use the BBHT algorithm in the way
we have described, we have to suppose that the number of elements of the
list is $N=2^{n}$ and we have to search for the indices of the elements
(instead for the elements themselves), because the domain of function $h$ in
Grover's algorithm is $\{0,\cdots ,N-1\}$. The details can be found in Ref.~%
\cite{DH99}. 
Dürr and H{ø}yer showed that the running time of the algorithm is $O(\sqrt{N}%
)$ and the probability of finding the minimum is at least 1/2. Their
analysis is valid when all elements of the list are distinct.

Baritompa \textit{et.~al} (BBW)~\cite{BBW05} used the DH algorithm to
propose a generic structure of a quantum global optimization algorithm for a
discrete function $f:\{0,\cdots ,N-1\}\rightarrow L$, where $L$ is a list of 
$N$ numbers, such that $f(x)$ is the $(x+1)$-th element in $L$. We use the
notation $h$ for the oracle function and $f$ for the function the minimum of
which we want to find. Note that in the DH algorithm, the BBHT algorithm is
used as a black box. BBW uses the BBHT algorithm explicitly, and can be
written in a pseudo-code as Algorithm~\ref{alg:GAS}. GAS is the basis for
the quantum algorithm used in this work.

\begin{algorithm}[!htb]

    \Begin{
         Generate $x_0$ uniformly in $\{0,\cdots,N-1\}$ and set $y_0=f(x_0)$.\;
         \For{ $i\,=1,2,\cdots$, until a termination condition is met} {
               Perform $r_{i}$ Grover's rotations marking points with image $\le y_{i-1}$. Denote outputs by $x$ and $y$\;
               \eIf{$y< y_{i-1}$}{
                          set $x_{i}=x$ and $y_{i}=y$\;
               }{
                          set $x_{i}=x_{i-1}$ and $y_{i}=y_{i-1}$\;
               }
         }
     }
\caption{Grover Adaptive Search (GAS)}
\label{alg:GAS}
\end{algorithm}

GAS reduces to the DH algorithm if: (1) the integer $r_{i}$ is chosen
uniformly at random in the range $0\leq r_{i}<m$, (2) $m$ is incremented as $%
m=\lambda m$ if $y\geq y_{i}$ and $m=1$ otherwise, where $\lambda =8/7,$ as
in the BBHT algorithm, and (3) the termination condition is that the
accumulated number of Grover's rotations is greater than $22.5\sqrt{N}%
+1.4\log ^{2}N$.

Baritompa \textit{et.~al} improved the prefactor that describes the running
time of the BBHT algorithm. They prove that for $\lambda =1.34$, the
expected number of oracle queries for the BBHT algorithm to find and verify
a marked element repeated $t$ times is at most $1.32\sqrt{N/|M|}$ \big(BBHT
uses the threshold $8\sqrt{N/|M|}$\big). They also provide a detailed proof
of the quadratic speedup of the DH algorithm when there are repeated
elements. They have proposed a new version of the minimization algorithm by
changing the method of choosing the number of Grover's rotations in each
round of the algorithm. Instead of selecting $r_{i}$ at random, they have
proposed a deterministic method in such way that the number of rotations
follows a pre-computed sequence of integers. The BBW version avoids to set $%
m=1$ each time the algorithm finds a new candidate for minimum. This is also
proposed in Ref.~\cite{KLPF08}, which shows that the number of measurement
reduces from $O(\log ^{2}N)$ in the DH algorithm to $O(\log N)$ in the
version that set $m=1$ only once at the beginning of the algorithm. The
running time of the BBW deterministic version is $2.46\sqrt{N},$ when there
are no repeated elements while the running time of the DH algorithm is $22.5%
\sqrt{N}$. There is no way to improve the scaling of those algorithms using
quantum computing, as been proved by Bennett \textit{et. al}~\cite{BBBV97}
and Zalka~\cite{Zal99}. Only the prefactor may be reduced.

\section[A New Method for Continuous Functions]{A New Method for Continuous
Functions}

\label{sec:3}

The new method proposed in this paper is a hybrid algorithm that employs a
classical optimization routine to find a local minimum and the GAS algorithm
to escape from that minimum towards another better candidate. We consider
continuous and differentiable real functions $f:D\rightarrow \mathbb{R}$
with $n$ variables, where the domain $D$ is a $n$-dimensional finite box. In
order to implement a computer program to find the global minimum point of a
continuous function, we discretize the function domain using intervals of
same length $\epsilon $ for all variables and we convert the domain points,
which form a $n$-dimensional array, into a one-dimensional array, generating
a list of $N$ points. After this discretization and conversion to
one-dimensional representation, the domain of $f$ can be taken as the set $%
\{0,\cdots ,N-1\}$. This process is mandatory because to perform Grover's
rotations we need an oracle function $h$ with domain $\{0,\cdots ,N-1\}$.
The value of $\epsilon $ depends on the structure of the function and on the
optimization problem. This parameter will be used in classical optimization
routines to characterize the precision of local minimum points, and at the
end to characterize the precision of the global minimum point. The general
structure of the new algorithm is given in Algorithm~\ref{alg:newmethod}.

\begin{algorithm}[!htb]
\Begin{
     Generate $x'$ uniformly at random in $\{0,\cdots,N-1\}$\;
     Use a classical optimization routine with input $x'$ to find a local minimum $x_0$ and set $y_0=f(x_0)$\;\label{line:3}
     Set $m=1$ and $\lambda=1.34$ (as suggested in BBW)\;
     \For{ $i\,=1,\,2,\cdots$, until a termination condition is met}{
          Define $M_i=\{x\in\{0,\cdots,N-1\}| f(x)<y_{i-1}\}$\;
          Choose $r_i$ uniformly at random in $\{0,\cdots,\lceil m-1\rceil\}$\;
          Apply $r_i$ Grover's rotations\;
          Perform a measurement. Let $x'\in\{0,\cdots,N-1\}$ be the output\;
          \eIf{ $x'\in M_i$}{
                    Use the classical optimization routine with input $x'$ to find a local minimum $x_i$\;\label{line:10}
                    Set $y_i=f(x_i)$\;
          }{
                    Set $x_i=x_{i-1}$, $y_i=f(x_i)$, and $m=\min\{\lambda\,m,\sqrt N\}$\;
          }
        }
     \Return last value of $x_i$\;
}
\caption{The New Method}
\label{alg:newmethod}
\end{algorithm}

The \textit{termination condition} determines the running time of the
algorithm. In the new method, the termination condition takes into account
the total number of Grover's rotations and the total number that the
objective function is evaluated when classical optimization routines are
employed. The weight of the classical objective function evaluation is
higher by a factor of $\sqrt{N}/\log N$ compared to the function evaluation
in each Grover's rotation. If $n_{1}$ is the number of Grover's rotations
and $n_{2}$ is the number of objective function evaluations in classical
optimization routines, then the termination condition is 
\begin{equation}\label{eq:term_cond}
n_{1}+\frac{\sqrt{N}}{\log^n N}\,n_{2}>2.46\sqrt{N},
\end{equation}%
where $n$ is the number of variables of the objective function.
In the worst case, the algorithm has running time $O(\sqrt{N})$, and in
cases for which classical algorithms are able to find the global
minimum without using the quantum part, the running time is $O(\log^n N)$. 

Algorithm~\ref{alg:newmethod} is a randomized algorithm. The success
probability cannot be calculated until the classical optimization
algorithm is specified. We are assuming that in the worst case the
classical algorithm will take $O(\log^n N)$ steps to find a local 
minimum after a initial point is given. The computational results 
presented in the next Section confirm that assumption.

\section[Computational Results]{Computational Results}

\label{sec:results}

Baritompa \textit{et.~al} compare their improved version of the minimization
algorithm (BBW) to the one of Dürr and H{ø}yer's (DH) by displaying
performance graphs, that depict the success probability of the algorithms in
terms of the (computational) effort. The effort is the number of objective
function evaluations before a new candidate for minimum point is found plus
the number of measurements. The number of measurements does not play an
important role for large $N$, because it scales logarithmically in terms of $%
N$. In the first part of this Section, we use the same technique to show
that the new method generates better results. This kind of analysis was
performed in Ref.~\cite{Hendrix:2000:0925-5001:143} for adaptive random
search, the structure of which is similar to the Grover adaptive search. The
details about the comparison between BBW and DH algorithms using performance
graphs can be obtained in Refs.~\cite{BBW05,LK10,LK11}.

\begin{figure}[!htb]
\begin{center}
\includegraphics[width=7.5cm]{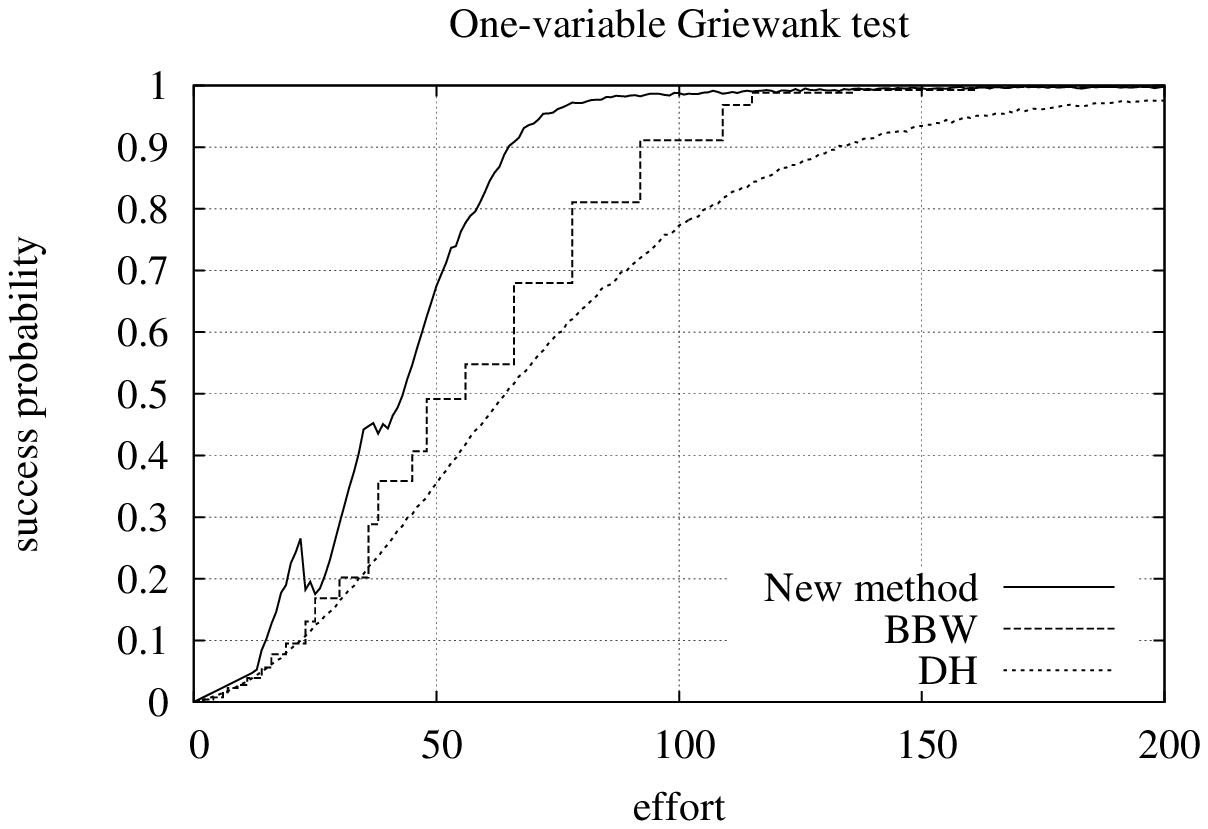}\quad %
\includegraphics[width=7.5cm]{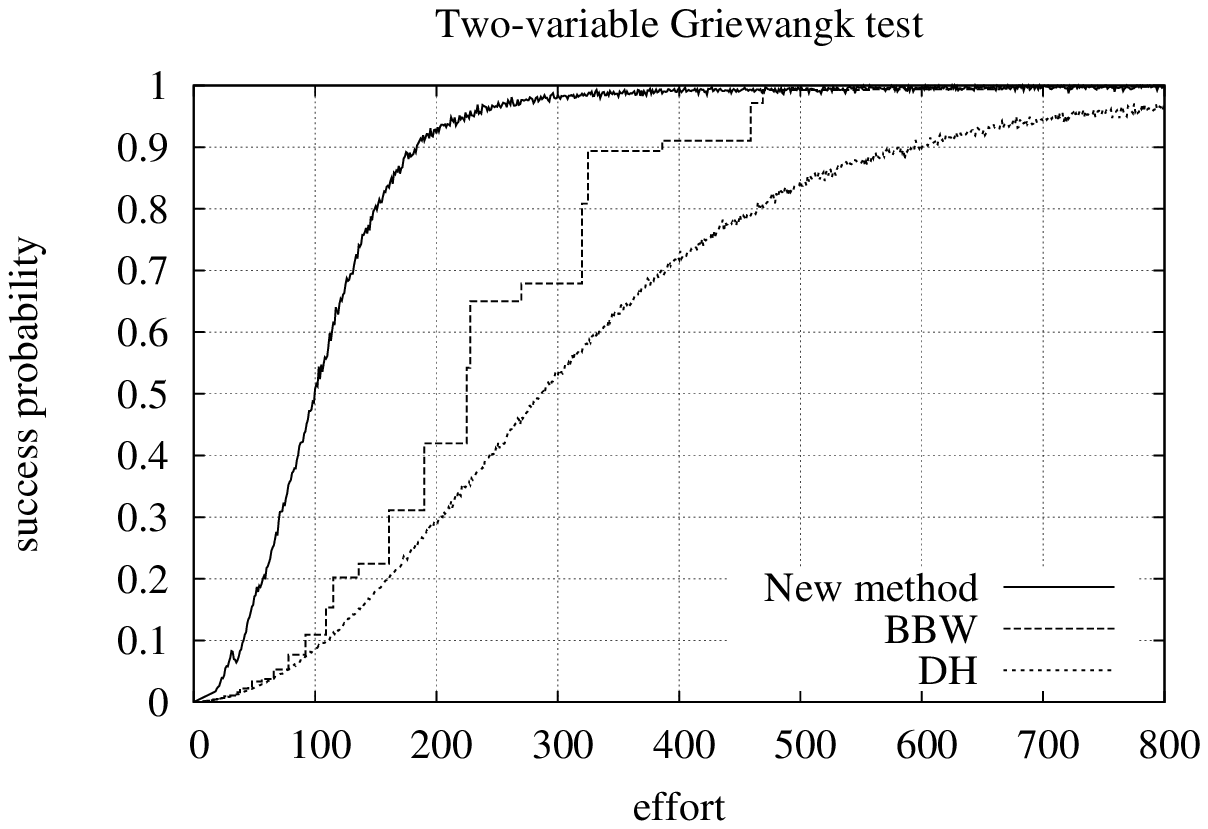}
\end{center}
\caption{Performance graphs comparing the new method with Baritompa \textit{%
et.~al}.~(BBW) and Dürr-Hoyer (DH) using one and two-variable Griewank test
functions.}
\label{fig_1and2vars}
\end{figure}

\begin{figure}[!htb]
\begin{center}
\includegraphics[width=7.5cm]{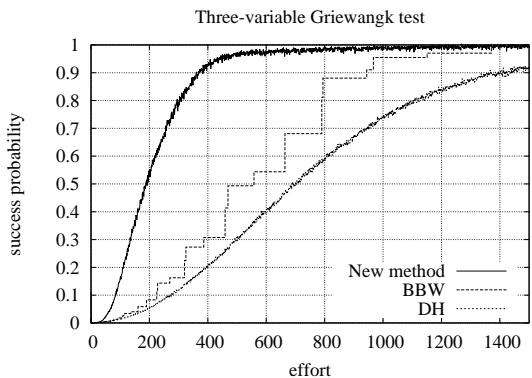}
\end{center}
\caption{Performance graphs comparing the new method with Baritompa \textit{%
et.~al}.~(BBW) and Dürr-Hoyer (DH) using three-variable Griewank test
function.}
\label{fig_3vars}
\end{figure}

Figs.~\ref{fig_1and2vars} and~\ref{fig_3vars} show the performance graphs
that compare the new method with the BBW and DH algorithms. We use one, two,
and three-variable Griewank test functions with domain $-40\leq
x_{0},x_{1},x_{2}\leq 40$, which are described in Appendix~\ref{appendix1}.
The classical routine used in Algorithm~\ref{alg:newmethod} (lines \ref%
{line:3} and \ref{line:10}) to find a local minimum is the BOBYQA routine~%
\cite{m8}. To generate those graphs, we create a sample with $N$ function
values taking $\{0,\ldots ,$ $N-1\}$ as the domain set (as described in Sec.~%
\ref{sec:3}). The value of $N$ for one-variable Griewank test is $N=2048$,
for two-variable is $N=256^{2}$, and for three-variable is $N=64^{3}$. We
average out this process 100,000 times for each graph. We take larger
parameter $\epsilon $ for the three-variable case, because to calculate the
average is time-consuming. All algorithms are implemented in the C language,
and it takes about half an hour in a 2.2 GHz Intel Core i7 processor to
generate the graphs of Fig.~\ref{fig_3vars}.

The new method is better than the previous methods in all those tests. Note
that it becomes even better when the number of variables increases, as can
be seen in Fig.~\ref{fig_3vars}. To have a success probability of 50\% for
the three-variable Griewank function, the effort of the new method is less
than 200 and the effort of the BBW and DH algorithms must be greater than
500. Similar curves are obtained when we use other classical optimization
routines, mainly for one and two-variable objective functions. For
three-variable functions, some optimization routines did not produce good
results. We discuss this problem later on, when we compare the efficiency of
different optimization routines.

To compare the new method with BBW and DH in more details, we employ the
most used performance test problems in global optimization, selecting test
functions with one to three variables (the list of test functions is in
Appendix~\ref{appendix1}). For our experiments we use the \textit{NLopt C
library}~\cite{nloptref}. NLopt is an open-source library under GNU license
for nonlinear optimization with a collection of classical optimization
routines. It can be used both for global and local optimization problems. In
this work, we use this library for the latter case. It is written in C and
has interface with many languages. The specific routines that we use in this
work are:

\begin{itemize}
\item LBFGS - This routine is based on variable-metric updates via Strang
recurrence of the low-storage BFGS routine~\cite{m1_a, m1_b}.

\item TNEWT - This routine is based on truncated Newton algorithms for
large-scale optimization~\cite{m2}.

\item MMA - This routine is based on a globally-convergent
method-of-moving-asymptotes algorithm for gradient-based local optimization,
including nonlinear inequality constraints~\cite{m3}.

\item COBYLA - This routine is a constrained optimization by linear
Approximations algorithm for derivative-free optimization with nonlinear
inequality and equality constraints~\cite{m4}.

\item NMEAD - This routine is a simplex method for function minimization
(Nelder-Mead simplex algorithm)~\cite{m5_a, m5_b}.

\item AGL - This routine is a globally convergent augmented Lagrangian
algorithm with general constraints and simple bounds~\cite{m7_a,m7_b}.

\item BOBYQA - This routine performs derivative-free bound-constrained
optimization using an iteratively constructed quadratic approximation for
the objective function~\cite{m8}.

\item CCSA - This routine is a conservative convex separable approximation
method, which is a variant of MMA algorithm~\cite{m3}.
\end{itemize}

In the experiments that compare the optimization routines, we create a
sample with $N$ function values taking $\{0,\ldots ,$ $N-1\}$ as the domain
set (as described in Sec.~\ref{sec:3}). The value of $N$ depends on the
number of variables: for functions with one variable, we take $N=2048$, for
2 variables, $N=2048^2$, and for 3 variables, $N=256^3$. With this sample,
we compute the running time (number of objective function evaluations) of
each algorithm until they find the correct global minimum. We are overriding
the termination conditions in all algorithms in order to check what is the
total number of function evaluations until the correct minimum is found. We
average out this process 100 times for each routine.

\noindent 
\begin{table}[!htb]
\begin{center}
{\footnotesize 
\begin{tabular}{|l|c|c|c|c|c|c|c|c|c|}
\hline
\diagbox{fcn.}{rout.} & LBFGS$^{*}$ & TNEWT$^{*}$ & MMA$^{*}$ & COBYLA$^{}$
& NMEAD$^{}$ & SBPLX$^{}$ & AGL$^{*}$ & BOBYQA$^{}$ & CCSAQ$^{*}$ \\ \hline
Neumaier & \textbf{\color{blue} 9.00 } & 17.52 & 33.00 & 11.00 & 193.5 & 
\textbf{\color{red} 236.2} & 34.00 & 11.00 & 32.67 \\ \hline
Griewank & 58.63 & 92.77 & 104.0 & 99.21 & \textbf{\color{red} 252.5} & 214.5
& 108.8 & \textbf{\color{blue} 52.61 } & 108.6 \\ \hline
Shekel & 12.00 & 15.00 & 32.67 & 11.81 & \textbf{\color{blue} 5.05 } & 11.37
& \textbf{\color{red} 33.66} & 11.00 & 33.00 \\ \hline
Rosenbrock & \textbf{\color{blue} 84.86 } & 110.5 & 153.8 & 129.7 & \textbf{%
\color{red} 220.4} & 207.9 & 159.1 & 101.7 & 165.5 \\ \hline
Michalewicz & 55.93 & 92.71 & 85.75 & 131.2 & \textbf{\color{red} 185.7} & 
153.6 & 89.46 & \textbf{\color{blue} 39.06 } & 91.82 \\ \hline
Dejong & \textbf{\color{blue} 9.00 } & 15.00 & 32.67 & 84.74 & 133.7 & 
\textbf{\color{red} 179.7} & 33.66 & 11.00 & 33.00 \\ \hline
Ackley & 90.39 & 106.3 & 129.6 & 136.6 & \textbf{\color{red} 212.8} & 210.3
& 141.6 & \textbf{\color{blue} 44.49 } & 145.4 \\ \hline
Schwefel & 9.00 & 12.00 & \textbf{\color{red} 66.26} & 21.95 & \textbf{%
\color{blue} 4.00 } & 25.15 & 34.00 & 19.57 & 33.00 \\ \hline
Rastrigin & 75.75 & 91.58 & \textbf{\color{blue} 33.00 } & 70.12 & \textbf{%
\color{red} 244.9} & 236.6 & 76.09 & 36.51 & 89.57 \\ \hline
Raydan & \textbf{\color{blue} 25.87 } & 43.97 & 33.00 & 55.94 & \textbf{%
\color{red} 191.8} & 178.6 & 34.00 & 28.89 & 36.73 \\ \hline
\end{tabular}
}
\end{center}
\caption{Average number of evaluations of one-variable testbed functions
using optimization routines specified in the first line. Smallest number of
evaluations in blue and largest in red.}
\label{table:onevar}
\end{table}

Table~\ref{table:onevar} shows the total number of function evaluations for
the one-variable case performed by the new method. The first line lists the
optimization routine used in Algorithm~\ref{alg:newmethod} (routines with
asterisk use function derivative) and the first column lists the name of the
test function. The expression and domain of the test functions are described
in Appendix~\ref{appendix1}. To highlight the best and worst routines
described in Table~\ref{table:onevar}, we highlight the smallest number of
evaluations in blue and largest in red for each test function. The results
show that the routine efficiency depends on the test function in general.
For one-variable functions, the LBFGS and BOBYQA routines have the best
performance while NMEAD has the worst, except for the Shekel and Schwefel
functions. The timings increase when we increase the number of sample values
(smaller $\epsilon $), but the table structure remains the same. On the
other hand, if we decrease too much the number of sample values, the
structure of the table changes significantly, and the information about the
best and worst routines is almost meaningless. From the timings in Table~\ref%
{table:onevar}, we cannot conclude that routines using function derivative
are better than routines that do not use derivative. The total simulation
time to produce the data in Table~\ref{table:onevar} is at order of some
minutes.

\begin{table}[!htb]
\begin{center}
{\footnotesize 
\begin{tabular}{|l|c|c|c|c|c|c|c|c|c|}
\hline
\diagbox{fcn.}{rout.} & LBFGS$^{*}$ & TNEWT$^{*}$ & MMA$^{*}$ & COBYLA$^{}$
& NMEAD$^{}$ & SBPLX$^{}$ & AGL$^{*}$ & BOBYQA$^{}$ & CCSAQ$^{*}$ \\ \hline
Neumaier & 1126 & 1671 & 1434 & \textbf{\color{blue} 22.00 } & 8943 & 
\textbf{\color{red} 9771} & 1204 & 122.8 & 1213 \\ \hline
Griewank & 1583 & 1441 & 1290 & 612.0 & \textbf{\color{red} 9285} & 9077 & 
1012 & \textbf{\color{blue} 412.2 } & 1036 \\ \hline
Shekel & 2643 & \textbf{\color{red} 2785} & 1575 & 22.00 & \textbf{%
\color{blue} 8.00 } & 22.01 & 1824 & 21.78 & 1933 \\ \hline
Rosenbrock & \textbf{\color{blue} 217.3 } & 770.2 & 1333 & 5762 & 7614 & 
\textbf{\color{red} 9707} & 1300 & 1086 & 1155 \\ \hline
Michalewicz & 1707 & 1881 & 1409 & 3887 & \textbf{\color{red} 10323} & 9841
& 1283 & \textbf{\color{blue} 390.4 } & 1258 \\ \hline
Dejong & 1534 & 1316 & 1379 & 22.00 & 6174 & \textbf{\color{red} 7356} & 1477
& \textbf{\color{blue} 21.78 } & 1433 \\ \hline
Ackley & 2385 & 2083 & 1487 & 1958 & \textbf{\color{red} 7611} & 7232 & 2001
& \textbf{\color{blue} 638.0 } & 2822 \\ \hline
Schwefel & \textbf{\color{red} 2101} & 1795 & 1422 & 356.7 & \textbf{%
\color{blue} 53.71 } & 1502 & 1306 & 70.32 & 1475 \\ \hline
Rastrigin & 1390 & 1219 & 1413 & 1360 & \textbf{\color{red} 7695} & 6180 & 
1265 & \textbf{\color{blue} 191.4 } & 1153 \\ \hline
Raydan & 1718 & 1726 & 1445 & 467.0 & 9430 & \textbf{\color{red} 10211} & 
1540 & \textbf{\color{blue} 324.7 } & 1651 \\ \hline
\end{tabular}
}
\end{center}
\caption{Average number of evaluations of two-variable testbed functions
using optimization routines specified in the first line. Smallest number of
evaluations in blue and largest in red.}
\label{table:twovar}
\end{table}

Table~\ref{table:twovar} shows the total number of function evaluations for
the two-variable case performed by the new method. The results show that the
best-performance routines depend heavily on the test function. The BOBYQA
routine has the best performance, while SBPLX has the worst performance. It
is remarkable that the NMEAD routine has the best performance for the Shekel
and Schwefel test functions and the worst performance for the Griewank,
Michalewicz, Ackley, and Rastrigin functions. This behavior also occurs with
one-variable functions, as we remarked earlier. All those six functions
are multimodal. The total simulation time to produce the
data in Table~\ref{table:twovar} is about one hour.

\begin{table}[!htb]
\begin{center}
{\footnotesize 
\begin{tabular}{|l|c|c|c|c|c|c|c|c|c|}
\hline
\diagbox{fcn.}{rout.} & LBFGS$^{*}$ & TNEWT$^{*}$ & MMA$^{*}$ & COBYLA$^{}$
& NMEAD$^{}$ & SBPLX$^{}$ & AGL$^{*}$ & BOBYQA$^{}$ & CCSAQ$^{*}$ \\ \hline
Neumaier & 1789 & 4467 & 3622 & 4147 & 6099 & \textbf{\color{red} 6805} & 
4080 & \textbf{\color{blue} 1390 } & 3215 \\ \hline
Griewank & 2932 & 2172 & 1500 & 1128 & 11392 & \textbf{\color{red} 12001} & 
3179 & \textbf{\color{blue} 711.4 } & 938.7 \\ \hline
Shekel & 4775 & \textbf{\color{red} 5137} & 1900 & 24.00 & \textbf{%
\color{blue} 12.00 } & 24.00 & 1300 & 24.00 & 1742 \\ \hline
Rosenbrock & \textbf{\color{blue} 1776 } & 1847 & 1898 & 3795 & 18074 & 
\textbf{\color{red} 18199} & 2592 & 4018 & 2531 \\ \hline
Michalewicz & 8051 & 4843 & -- & 4037 & \textbf{\color{red} 14779} & 9944 & 
2209 & 3873 & \textbf{\color{blue} 1468 } \\ \hline
Dejong & 2604 & 2685 & 1582 & 1005 & \textbf{\color{red} 9735} & 6608 & 1953
& \textbf{\color{blue} 23.76 } & 1338 \\ \hline
Ackley & -- & \textbf{\color{blue} 918.3 } & 1674 & 2944 & \textbf{%
\color{red} 16455} & 15680 & 2390 & 1098 & 1249 \\ \hline
Schwefel & 6988 & \textbf{\color{red} 7549} & 557.0 & \textbf{\color{blue}
141.7 } & 1100 & 209.3 & 883.2 & 172.4 & 1619 \\ \hline
Rastrigin & 2578 & \textbf{\color{blue} 555.7 } & 1779 & 1638 & \textbf{%
\color{red} 11386} & 9408 & 1365 & 925.1 & 1180 \\ \hline
Raydan & 2826 & 2908 & 1816 & \textbf{\color{blue} 726.9 } & \textbf{%
\color{red} 17283} & 16469 & 2081 & 780.1 & 1628 \\ \hline
\end{tabular}%
}
\end{center}
\caption{Average number of evaluations of three-variable testbed functions
using optimization routines specified in the first line. Smallest number of
evaluations in blue and largest in red.}
\label{table:threevar}
\end{table}

Table~\ref{table:threevar} shows the total number of function evaluations
for the three-variable case performed by the new method. The results show
again that the best-performance routines depend heavily on the test
function. The BOBYQA routine has a small advantage while NMEAD and SBPLX
have the worst performance. Notice that we are using $N=256^3$ in the
discretization procedure, which means that we use 256 points in each axis.
This number is rather small and represents the continuous functions in a
gross manner. We do not use larger number of function values, because the
total simulation time, which includes the averages, is long. Using such
small number of function values, the best method seems to be BOBYQA and the
worst NMEAD, similar to what happens with two-variable functions. It is also
similar to what happens with one-variable functions, except that LBFGS is
not as efficient as in one-variable case. The total simulation time to
produce the data in Table~\ref{table:threevar} is about five hours.

\begin{table}[!htb]
\begin{center}
{\small 
\begin{tabular}{|l|l|l|l|l|l|l|l|l|l|}
\hline
\textbf{\# variables} & \multicolumn{3}{|c|}{\textbf{one variable}} & 
\multicolumn{3}{|c|}{\textbf{two variables}} & \multicolumn{3}{|c|}{\textbf{%
three variables}} \\ \hline
\diagbox{fcn.}{meth.} & NEW & BBW & DH & NEW & BBW & DH & NEW & BBW & DH \\ 
\hline
Neumaier & \textbf{\color{blue} 9.00 } & 97.50 & \textbf{\color{red} 112.4}
& \textbf{\color{blue} 22.00 } & \textbf{\color{red} 960.8} & 694.0 & 1390 & 
\textbf{\color{blue} 371.3 } & \textbf{\color{red} 13660} \\ \hline
Griewank & \textbf{\color{blue} 52.61 } & 70.21 & \textbf{\color{red} 88.44}
& \textbf{\color{blue} 412.2 } & 865.1 & \textbf{\color{red} 1119} & \textbf{%
\color{blue} 711.4 } & 884.4 & \textbf{\color{red} 2643} \\ \hline
Shekel & \textbf{\color{blue} 5.05 } & 99.04 & \textbf{\color{red} 113.6} & 
\textbf{\color{blue} 8.00 } & 944.0 & \textbf{\color{red} 4588} & \textbf{%
\color{blue} 12.00 } & 528.0 & \textbf{\color{red} 18143} \\ \hline
Rosenbrock & \textbf{\color{blue} 84.86 } & 87.68 & \textbf{\color{red} 124.3%
} & \textbf{\color{blue} 217.3 } & 852.8 & \textbf{\color{red} 11017} & 
\textbf{\color{blue} 1776 } & -- & \textbf{\color{red} 4580} \\ \hline
Michalewicz & \textbf{\color{blue} 39.06 } & 90.13 & \textbf{\color{red}
107.7} & \textbf{\color{blue} 390.4 } & 739.3 & \textbf{\color{red} 2794} & 
1468 & \textbf{\color{blue} 685.5 } & \textbf{\color{red} 17334} \\ \hline
Dejong & \textbf{\color{blue} 9.00 } & 75.00 & \textbf{\color{red} 97.50} & 
\textbf{\color{blue} 21.78 } & 826.6 & \textbf{\color{red} 3056} & \textbf{%
\color{blue} 23.76 } & 826.9 & \textbf{\color{red} 5886} \\ \hline
Ackley & \textbf{\color{blue} 44.49 } & 102.0 & \textbf{\color{red} 114.4} & 
\textbf{\color{blue} 638.0 } & 805.7 & \textbf{\color{red} 875.0} & 918.3 & 
\textbf{\color{blue} 617.0 } & \textbf{\color{red} 4260} \\ \hline
Schwefel & \textbf{\color{blue} 4.00 } & 88.38 & \textbf{\color{red} 124.2}
& \textbf{\color{blue} 53.71 } & 871.1 & \textbf{\color{red} 3675} & \textbf{%
\color{blue} 141.7 } & 685.5 & \textbf{\color{red} 13475} \\ \hline
Rastrigin & \textbf{\color{blue} 33.00 } & 74.48 & \textbf{\color{red} 99.45}
& \textbf{\color{blue} 191.4 } & 792.2 & \textbf{\color{red} 3510} & \textbf{%
\color{blue} 555.7 } & 967.0 & \textbf{\color{red} 1010} \\ \hline
Raydan & \textbf{\color{blue} 25.87 } & 93.24 & \textbf{\color{red} 115.2} & 
\textbf{\color{blue} 324.7 } & 609.6 & \textbf{\color{red} 10063} & \textbf{%
\color{blue} 726.9 } & 774.0 & \textbf{\color{red} 4948} \\ \hline
\end{tabular}
}
\end{center}
\caption{Total number of objective function evaluations for the new method,
the BBW and DH algorithms using one, two, and three-variable test functions.}
\label{table:comparison}
\end{table}

Table~\ref{table:comparison} shows the total number of objective function
evaluations until the global minimum is found for the three algorithms used
in this work. In this experiment, we use the same discretization parameters
used in the previous tables and we override again the termination condition
of the algorithms, that is, the algorithms run until the correct global
minimum is found. When we use the correct termination condition for each
algorithm, the number of objective function evaluations is close to the ones
showed in the Table~\ref{table:comparison}, but in a fraction of cases
(smaller than 50\%) we do not find the correct global minimum point. In this
experiment, the new method uses the most efficient routine for each test
function, information that is obtained from Tables~\ref{table:onevar}, 
\ref{table:twovar}, and~\ref{table:threevar}. For one and two-variable test
functions, the advantage of the new method is remarkable. The advantage in
the three-variable case is not so impressive, in contradiction to what we
have concluded earlier, when we analyzed the performance graphs in Figs.~\ref%
{fig_1and2vars} and~\ref{fig_3vars}. For the Neumaier, Michalewicz, and
Ackley functions, the BBW algorithm is better than the new method. This
seems to be a consequence of the small number of functions values in the
discretization procedure ($N=256^{3}$).

\begin{table}[!htb]
\begin{center}
{\small 
\begin{tabular}{|l|c|c|c|} \hline
\diagbox{fcn.}{\# vars.} & one-variable & two-variable & three-variable \\ \hline
Neumaier & 1.00 & 1.00 & 1.00 \\\hline 
Griewank & 0.87 & 1.00 & 1.00 \\\hline 
Shekel & 0.87 & 0.96 & 1.00 \\\hline 
Rosenbrock & 0.99 & 1.00 & 0.99 \\\hline 
Michalewicz & 1.00 & 1.00 & 0.99 \\\hline 
Dejong & 0.97 & 0.99 & 1.00 \\\hline 
Ackley & 1.00 & 1.00 & 1.00 \\\hline 
Schwefel & 0.98 & 1.00 & 1.00 \\\hline 
Rastrigin & 1.00 & 1.00 & 1.00 \\\hline 
Raydan & 1.00 & 1.00 & 1.00 \\\hline\hline 
{AVERAGE} & 0.96 & 0.99 & 0.99 \\\hline  
 \end{tabular}
}
\end{center}
\caption{Success probability of Algorithm~\ref{alg:newmethod} using the
best classical optimization routine.}
\label{table:probability_mail}
\end{table}

Table~\ref{table:probability_mail} shows the success probability of 
Algorithm~\ref{alg:newmethod} when we select the best optimization routine.
The success probability is calculated in the following way: We run
Algorithm~\ref{alg:newmethod} many times 
with the termination condition given by Eq.~(\ref{eq:term_cond})
and we count the number of times that the algorithm
finds the correct global minimum. The success probability is the
success rate. Table~\ref{table:probability_mail} is built using
the data described in Appendix~\ref{appendix3}, which shows the
success probability of Algorithm~\ref{alg:newmethod} for all
optimization routines. Notice that the success probabilities 
improve when we increase the number of variables. This shows that
the termination condition given by Eq.~(\ref{eq:term_cond}) is
sound. If the success probability for one variable (0.94) is not
high enough for practical purposes, one can rerun the algorithm
many times in order to improve this value. 
From Appendix~\ref{appendix3},
we conclude that the BOBYQA routine is the best one on average for the 
testbed functions with one to three variables,
while the SBPLX routine is the worst on average for one and two-variable functions and
the NMEAD routine is the worst on average for three-variable functions. 
As an exception, the BOBYQA routine has a bad performance for the one-variable
Rosenbrock function. Similar conclusions
were drawn from Tables~\ref{table:onevar}, \ref{table:twovar}, and~\ref{table:threevar}.
Those coincidences were expected since the experiments are correlated.

\section{Conclusions}

\label{sec:conclusion}

This paper proposed a new method for continuous global optimization
problems, using GAS and classical routines to find efficiently a local
minimum. Our numerical simulations show that DH, BBW, and the new method
have very different asymptotic behavior, where the new method presented a
better performance.

\section*{Acknowledgments}

The authors would like to thank FAPESP and CNPq for their financial support.
R.P. would like to thank prof.~Benjamín Barán for useful suggestions.

\newpage

\appendix

\section{Test functions}
\label{appendix1}

{\noindent{\large Neumaier}} 
\begin{equation*}
f(x_0,\ldots,x_{n-1})=\sum_{i=0}^{n-1} (x_i - 1)^2 - \sum_{i=1}^{n-1} x_i
x_{i-1}, \ 0 \leq x_i\leq 4
\end{equation*}

{\noindent{\large Griewank}} 
\begin{equation*}
f(x_0,\ldots,x_{n-1})=\frac{1}{4000} \sum_{i = 0} ^ {n-1} x_i^2 -
\prod_{i=0}^{n-1} \cos \left( \frac{x_i}{\sqrt{i+1}} \right) + 1, \ -40 \leq
x_i\leq 40
\end{equation*}

{\noindent{\large Shekel}} 
\begin{equation*}
f(x_0,\ldots,x_{n-1}) = \sum_{i = 0}^{m-1} \frac{1}{c_{i} + \sum_{j =
0}^{n-1} (x_{j} - a_{ji})^2 }, \ -1 \leq x_i\leq 1
\end{equation*}

{\noindent{\large Rosenbrock}} 
\begin{equation*}
f(x_0,\ldots,x_{n-1}) = \sum_{i=0}^{n-2} (1-x_i)^2+ 100 (x_{i+1} - x_i^2 )^2
, \ -30 \leq x_i\leq 30
\end{equation*}

{\noindent{\large Michalewicz }} 
\begin{equation*}
f(x_0,\ldots,x_{n-1}) =-\sum_{i=0}^{n-1} \sin(x_i) \sin^{2m}\left(\frac{ i
x_i^2}{\pi}\right), \ 0 \leq x_i\leq 10
\end{equation*}

{\noindent{\large Dejong}} 
\begin{equation*}
f(x_0,\ldots,x_{n-1}) =\sum_{i=0}^{n-1}x_i^2, \ -5.12 \leq x_i\leq 5.12
\end{equation*}

{\noindent{\large Ackley}} 
\begin{equation*}
f(x_0,\ldots,x_{n-1}) = -20 \exp\left(-\frac{1}{5} \sqrt{\frac{1}{n}%
\sum_{i=0}^{n-1} x_i^2} \; \right) -\exp\left(\frac{1}{n} \sum_{i=0}^{n-1}
\cos( 2 \pi x_i) \right) + 20 +\exp(1), \ -15 \leq x_i\leq 20
\end{equation*}

{\noindent{\large Schwefel}} 
\begin{equation*}
f(x_0,\ldots,x_{n-1}) = -\sum_{i=0}^{n-1} x_i \sin\left(\sqrt{|x_i|}\right),
\ -20 \leq x_i\leq 20
\end{equation*}

{\noindent{\large Rastrigin}} 
\begin{equation*}
f(x_0,\ldots,x_{n-1}) = \sum_{i=0}^{n-1} \left( x_i^2 - 10\cos(2 \pi x_i) +
10 \right), \ -5.12 \leq x_i\leq 5.12
\end{equation*}

{\noindent{\large Raydan}} 
\begin{equation*}
f(x_0,\ldots,x_{n-1}) = -\sum_{i=0}^{n-1} \frac{(i+1)}{10}\left( \exp(x_i) -
x_i\right), \ -5.12 \leq x_i\leq 5.12
\end{equation*}

\newpage

\section{Standard Deviation}
\label{appendix2}

This Appendix shows the tables of the standard deviation of the number of
evaluations for one, two, and three-variable test functions associated with
Tables~\ref{table:onevar}, \ref{table:twovar}, and~\ref{table:threevar},
respectively. In all tables, the smallest standard deviations are depicted
in blue and largest in red.

\noindent 
\begin{table}[!htb]
{\footnotesize 
\begin{tabular}{|l|c|c|c|c|c|c|c|c|c|}
\hline
\diagbox{fcn.}{rout.} & LBFGS$^{*}$ & TNEWT$^{*}$ & MMA$^{*}$ & COBYLA$^{}$
& NMEAD$^{}$ & SBPLX$^{}$ & AGL$^{*}$ & BOBYQA$^{}$ & CCSAQ$^{*}$ \\ \hline
Neumaier & \textbf{\color{blue} 0.00 } & 3.55 & 0.00 & 0.00 & 150.4 & 
\textbf{\color{red} 210.3} & 0.00 & 0.00 & 3.27 \\ \hline
Griewank & 22.05 & 38.36 & 40.63 & 90.67 & \textbf{\color{red} 212.7} & 181.5
& 37.47 & \textbf{\color{blue} 21.55 } & 41.86 \\ \hline
Shekel & \textbf{\color{blue} 0.00 } & 0.00 & 3.27 & \textbf{\color{red} 9.28%
} & 2.82 & 2.60 & 3.37 & 0.00 & 0.00 \\ \hline
Rosenbrock & \textbf{\color{blue} 35.08 } & 45.33 & 52.07 & 96.28 & \textbf{%
\color{red} 190.4} & 146.1 & 55.62 & 52.44 & 50.81 \\ \hline
Michalewicz & 29.42 & 39.70 & 35.52 & 107.4 & \textbf{\color{red} 136.3} & 
103.4 & 36.53 & \textbf{\color{blue} 23.93 } & 38.09 \\ \hline
Dejong & \textbf{\color{blue} 0.00 } & 0.00 & 3.27 & 62.03 & 128.4 & \textbf{%
\color{red} 139.2} & 3.37 & 0.00 & 0.00 \\ \hline
Ackley & \textbf{\color{blue} 48.83 } & 54.81 & 54.75 & 91.30 & \textbf{%
\color{red} 182.4} & 171.6 & 54.63 & 50.23 & 52.63 \\ \hline
Schwefel & \textbf{\color{blue} 0.00 } & 0.00 & \textbf{\color{red} 26.05} & 
13.99 & 0.00 & 24.04 & 0.00 & 8.44 & 0.00 \\ \hline
Rastrigin & 28.23 & 32.11 & \textbf{\color{blue} 0.00 } & 40.89 & \textbf{%
\color{red} 216.2} & 207.0 & 34.35 & 17.76 & 34.00 \\ \hline
Raydan & 5.88 & 21.42 & \textbf{\color{blue} 0.00 } & 27.22 & \textbf{%
\color{red} 177.9} & 143.3 & 0.00 & 34.48 & 15.96 \\ \hline
\end{tabular}
}
\caption{Standard deviation for one-variable test functions.}
\label{table:onevar_sd}
\end{table}

\noindent 
\begin{table}[!htb]
{\footnotesize 
\begin{tabular}{|l|c|c|c|c|c|c|c|c|c|}
\hline
\diagbox{fcn.}{rout.} & LBFGS$^{*}$ & TNEWT$^{*}$ & MMA$^{*}$ & COBYLA$^{}$
& NMEAD$^{}$ & SBPLX$^{}$ & AGL$^{*}$ & BOBYQA$^{}$ & CCSAQ$^{*}$ \\ \hline
Neumaier & 1240 & 3163 & 849.5 & \textbf{\color{blue} 0.00 } & \textbf{%
\color{red} 9011} & 8935 & 712.7 & 452.5 & 602.9 \\ \hline
Griewank & 2757 & 1394 & 633.2 & 475.0 & 9369 & \textbf{\color{red} 10116} & 
531.5 & \textbf{\color{blue} 271.4 } & 412.8 \\ \hline
Shekel & \textbf{\color{red} 2736} & 2692 & 1076 & \textbf{\color{blue} 0.00 
} & 0.00 & 3.17 & 1305 & 2.18 & 1531 \\ \hline
Rosenbrock & \textbf{\color{blue} 307.4 } & 751.4 & 1258 & 6865 & \textbf{%
\color{red} 8471} & 7946 & 883.1 & 1713 & 714.2 \\ \hline
Michalewicz & 1817 & 2190 & 900.0 & 3552 & \textbf{\color{red} 9537} & 9376
& 670.2 & \textbf{\color{blue} 615.7 } & 658.8 \\ \hline
Dejong & 2625 & 917.5 & 909.1 & \textbf{\color{blue} 0.00 } & 7359 & \textbf{%
\color{red} 7635} & 1042 & 2.18 & 817.4 \\ \hline
Ackley & 1693 & 1778 & 1017 & 4867 & \textbf{\color{red} 7990} & 6736 & 1177
& \textbf{\color{blue} 918.7 } & 2156 \\ \hline
Schwefel & 4390 & 2498 & 1041 & 2038 & \textbf{\color{blue} 23.94 } & 
\textbf{\color{red} 4840} & 1030 & 182.6 & 851.1 \\ \hline
Rastrigin & 1006 & 825.5 & 856.9 & 1567 & \textbf{\color{red} 8918} & 6952 & 
922.0 & \textbf{\color{blue} 110.3 } & 679.0 \\ \hline
Raydan & 2823 & 2181 & 791.9 & \textbf{\color{blue} 493.0 } & \textbf{%
\color{red} 9429} & 8425 & 1068 & 535.3 & 1336 \\ \hline
\end{tabular}
}
\caption{Standard deviation for two-variable test functions.}
\label{table:twovar_sd}
\end{table}

\noindent 
\begin{table}[!htb]
{\footnotesize 
\begin{tabular}{|l|c|c|c|c|c|c|c|c|c|}
\hline
\diagbox{fcn.}{rout.} & LBFGS$^{*}$ & TNEWT$^{*}$ & MMA$^{*}$ & COBYLA$^{}$
& NMEAD$^{}$ & SBPLX$^{}$ & AGL$^{*}$ & BOBYQA$^{}$ & CCSAQ$^{*}$ \\ \hline
Neumaier & \textbf{\color{blue} 1421 } & 5056 & 3018 & 2847 & \textbf{%
\color{red} 8923} & 5652 & 3628 & 2173 & 2425 \\ \hline
Griewank & 3866 & 2720 & 801.1 & 1460 & \textbf{\color{red} 15369} & 15037 & 
4247 & 699.2 & \textbf{\color{blue} 695.4 } \\ \hline
Shekel & 6247 & \textbf{\color{red} 7481} & 1599 & \textbf{\color{blue} 0.00 
} & 0.00 & 0.00 & 843.6 & 0.00 & 1154 \\ \hline
Rosenbrock & 2196 & \textbf{\color{blue} 1306 } & 1390 & 3521 & \textbf{%
\color{red} 17103} & 15874 & 1806 & 3026 & 1931 \\ \hline
Michalewicz & 5928 & 2770 & -- & 4859 & \textbf{\color{red} 13681} & 8277 & 
2209 & 4610 & \textbf{\color{blue} 1468} \\ \hline
Dejong & 5084 & 3581 & 1082 & 843.7 & \textbf{\color{red} 14373} & 9545 & 
1241 & \textbf{\color{blue} 2.38 } & 857.6 \\ \hline
Ackley & -- & 1299 & 1674 & 3644 & \textbf{\color{red} 15637} & 14966 & 1797
& 1654 & \textbf{\color{blue} 1026} \\ \hline
Schwefel & 8141 & \textbf{\color{red} 9863} & 557.0 & \textbf{\color{blue}
343.5 } & 2164 & 645.1 & 550.4 & 870.5 & 1155 \\ \hline
Rastrigin & 2024 & 785.8 & 1266 & 1727 & 12457 & \textbf{\color{red} 12655}
& 868.4 & \textbf{\color{blue} 636.2 } & 798.6 \\ \hline
Raydan & 4042 & 3778 & 1325 & \textbf{\color{blue} 778.6 } & 15591 & \textbf{%
\color{red} 16436} & 1393 & 877.8 & 1058 \\ \hline
\end{tabular}
}
\caption{Standard deviation for three-variable test functions.}
\end{table}

\newpage

\section{Success Probability}
\label{appendix3}

This Appendix shows the tables of success probability of 
Algorithm~\ref{alg:newmethod} with the termination condition 
given by Eq.~(\ref{eq:term_cond})
for one, two, and three-variable test functions 
using the classical optimization routines. 
In all tables, the largest probability are depicted
in blue and lowest in red. We have performed an average over
100 rounds for each table.

\noindent 
\begin{table}[!htb]
{\footnotesize 
\begin{tabular}{|l|c|c|c|c|c|c|c|c|c|} \hline
\diagbox{fcn.}{rout.}
 & LBFGS$^{*}$ & TNEWT$^{*}$ & MMA$^{*}$ & COBYLA$^{}$ & NMEAD$^{}$ & SBPLX$^{}$ & AGL$^{*}$ & BOBYQA$^{}$ & CCSAQ$^{*}$ \\ \hline
Neumaier & {\bf\color{blue} 1.00} & {\bf\color{blue} 1.00} & {\bf\color{blue} 1.00} & 0.37 & {\bf\color{red} 0.10 } & 0.12 & {\bf\color{blue} 1.00} & {\bf\color{blue} 1.00} & {\bf\color{blue} 1.00} \\\hline 
Griewank & 0.85 & 0.59 & 0.42 & 0.83 & 0.06 & {\bf\color{red} 0.04 } & 0.42 & {\bf\color{blue} 0.87} & 0.40 \\\hline 
Shekel & {\bf\color{blue} 1.00}  & {\bf\color{blue} 1.00} & {\bf\color{blue} 1.00} & {\bf\color{blue} 1.00} & {\bf\color{blue} 1.00} & {\bf\color{blue} 1.00} & {\bf\color{blue} 1.00} & {\bf\color{blue} 1.00} & {\bf\color{blue} 1.00} \\\hline 
Rosenbrock & {\bf\color{blue} 0.87} & 0.50 & 0.25 & 0.44 & 0.13 & {\bf\color{red} 0.10 } & 0.25 & 0.27 & 0.15 \\\hline 
Michalewicz & 0.82 & 0.46 & 0.68 & 0.21 & {\bf\color{red} 0.08 } & 0.08 & 0.68 & {\bf\color{blue} 0.99} & 0.58 \\\hline 
Dejong & {\bf\color{blue} 1.00} & {\bf\color{blue} 1.00} & {\bf\color{blue} 1.00} & 0.63 & 0.22 & {\bf\color{red} 0.16 } & {\bf\color{blue} 1.00} & {\bf\color{blue} 1.00} & {\bf\color{blue} 1.00} \\\hline 
Ackley & 0.57 & 0.44 & 0.32 & 0.23 & 0.13 & {\bf\color{red} 0.10 } & 0.32 & {\bf\color{blue} 0.97} & 0.34 \\\hline 
Schwefel & 0.91 & 0.84 & {\bf\color{red} 0.79 } & {\bf\color{blue} 1.00} & {\bf\color{blue} 1.00} & {\bf\color{blue} 1.00} & 0.79 & {\bf\color{blue} 1.00} & 0.88 \\\hline 
Rastrigin & 0.52 & 0.53 & 0.76 & 0.51 & 0.03 & {\bf\color{red} 0.02 } & 0.76 & {\bf\color{blue} 0.98} & 0.59 \\\hline 
Raydan & {\bf\color{blue} 1.00} & {\bf\color{blue} 1.00} & {\bf\color{blue} 1.00} & 0.87 & 0.11 & {\bf\color{red} 0.08 } & {\bf\color{blue} 1.00} & {\bf\color{blue} 1.00} & {\bf\color{blue} 1.00} \\\hline\hline 
{AVERAGE} & 0.854 & 0.736 & 0.722 & 0.609 & 0.286 & {\bf\color{red} 0.27 } & 0.722 & {\bf\color{blue} 0.908} & 0.694 \\\hline 
\end{tabular}
}
\caption{Success probability for one-variable test functions.}
\label{table:onevar_sd}
\end{table}

\noindent 
\begin{table}[!htb]
{\footnotesize 
\begin{tabular}{|l|c|c|c|c|c|c|c|c|c|} \hline
\diagbox{fcn.}{rout.}
 & LBFGS$^{*}$ & TNEWT$^{*}$ & MMA$^{*}$ & COBYLA$^{}$ & NMEAD$^{}$ & SBPLX$^{}$ & AGL$^{*}$ & BOBYQA$^{}$ & CCSAQ$^{*}$ \\ \hline
Neumaier & 0.96 & 0.98 & 0.92 & 0.86 & 0.37 & {\bf\color{red} 0.28 } & 0.93 & {\bf\color{blue} 1.00} & 0.98 \\\hline 
Griewank & {\bf\color{blue} 1.00} & 0.94 & 0.81 & {\bf\color{blue} 1.00} & {\bf\color{red} 0.39 } & 0.39 & 0.84 & {\bf\color{blue} 1.00} & 0.86 \\\hline 
Shekel & 0.95 & 0.94 & 0.74 & {\bf\color{blue} 1.00} & {\bf\color{blue} 1.00} & {\bf\color{blue} 1.00} & {\bf\color{red} 0.70 } & {\bf\color{blue} 1.00} & 0.77 \\\hline 
Rosenbrock & 0.90 & 0.93 & 0.95 & 0.65 & 0.35 & {\bf\color{red} 0.30 } & {\bf\color{blue} 0.96} & 0.95 & 0.95 \\\hline 
Michalewicz & 0.96 & 0.97 & 0.87 & 0.98 & 0.42 & {\bf\color{red} 0.36 } & 0.86 & {\bf\color{blue} 1.00} & 0.90 \\\hline 
Dejong & 0.99 & 0.97 & 0.93 & {\bf\color{blue} 1.00} & 0.72 & {\bf\color{red} 0.52 } & 0.93 & {\bf\color{blue} 1.00} & 0.93 \\\hline 
Ackley & 0.71 & 0.46 & 0.68 & {\bf\color{blue} 0.99} & {\bf\color{red} 0.33 } & 0.37 & 0.79 & {\bf\color{blue} 0.99} & 0.81 \\\hline 
Schwefel & 0.97 & 0.95 & 0.85 & {\bf\color{blue} 1.00} & {\bf\color{blue} 1.00} & 0.94 & 0.87 & {\bf\color{blue} 1.00} & {\bf\color{red} 0.84 } \\\hline 
Rastrigin & 0.94 & 0.78 & 0.88 & {\bf\color{blue} 1.00} & 0.51 & {\bf\color{red} 0.41 } & 0.86 & {\bf\color{blue} 1.00} & 0.91 \\\hline 
Raydan & 0.99 & 0.99 & 0.94 & {\bf\color{blue} 1.00} & 0.34 & {\bf\color{red} 0.29 } & 0.90 & {\bf\color{blue} 1.00} & 0.94 \\\hline\hline 
{AVERAGE} & 0.937 & 0.891 & 0.857 & 0.948 & 0.543 & {\bf\color{red} 0.486 } & 0.864 & {\bf\color{blue} 0.994} & 0.889 \\\hline 
\end{tabular}
}
\caption{Success probability for two-variable test functions.}
\label{table:twovar_sd}
\end{table}

\noindent 
\begin{table}[!htb]
{\footnotesize 
\begin{tabular}{|l|c|c|c|c|c|c|c|c|c|} \hline
\diagbox{fcn.}{rout.}
 & LBFGS$^{*}$ & TNEWT$^{*}$ & MMA$^{*}$ & COBYLA$^{}$ & NMEAD$^{}$ & SBPLX$^{}$ & AGL$^{*}$ & BOBYQA$^{}$ & CCSAQ$^{*}$ \\ \hline
Neumaier & 0.94 & {\bf\color{red} 0.80 } & 0.88 & {\bf\color{blue} 1.00} & 0.96 & 0.91 & 0.89 & 0.98 & 0.90 \\\hline 
Griewank & 0.95 & 0.96 & 0.94 & {\bf\color{blue} 1.00} & {\bf\color{red} 0.57 } & 0.60 & 0.92 & {\bf\color{blue} 1.00} & 0.95 \\\hline 
Shekel & 0.95 & 0.92 & 0.93 & {\bf\color{blue} 1.00} & {\bf\color{blue} 1.00} & {\bf\color{blue} 1.00} & 0.92 & {\bf\color{blue} 1.00} & {\bf\color{red} 0.91 } \\\hline 
Rosenbrock & 0.69 & 0.75 & 0.85 & 0.72 & 0.43 & {\bf\color{red} 0.40 } & 0.89 & {\bf\color{blue} 0.99} & 0.88 \\\hline 
Michalewicz & 0.79 & 0.86 & 0.78 & 0.77 & {\bf\color{red} 0.59 } & 0.68 & 0.77 & {\bf\color{blue} 0.99} & 0.80 \\\hline 
Dejong & 0.99 & 0.99 & 0.97 & {\bf\color{blue} 1.00} & {\bf\color{red} 0.84 } & 0.86 & 0.99 & {\bf\color{blue} 1.00} & 0.92 \\\hline 
Ackley & 0.77 & 0.68 & 0.89 & 0.96 & {\bf\color{red} 0.46 } & 0.70 & 0.89 & {\bf\color{blue} 1.00} & 0.83 \\\hline 
Schwefel & 0.92 & 0.88 & 0.81 & 0.99 & {\bf\color{blue} 1.00} & 0.99 & {\bf\color{red} 0.73 } & {\bf\color{blue} 1.00} & 0.84 \\\hline 
Rastrigin & 0.78 & 0.87 & 0.95 & {\bf\color{blue} 1.00} & {\bf\color{red} 0.59 } & 0.64 & 0.94 & {\bf\color{blue} 1.00} & 0.96 \\\hline 
Raydan & 0.94 & 0.94 & 0.96 & 0.99 & {\bf\color{red} 0.39 } & 0.47 & 0.92 & {\bf\color{blue} 1.00} & 0.93 \\\hline\hline 
{AVERAGE} & 0.872 & 0.865 & 0.896 & 0.943 & {\bf\color{red} 0.683 } & 0.725 & 0.886 & {\bf\color{blue} 0.996} & 0.892 \\\hline 
\end{tabular}
}
\caption{Success probability for three-variable test functions.}
\end{table}

\end{document}